\documentclass[a4paper, 12pt]{article}

\usepackage[sort&compress]{natbib}
\bibpunct{(}{)}{;}{a}{}{,} 

\usepackage{amsthm, amsmath, amssymb, mathrsfs, multirow, url, subfigure}
\usepackage{graphicx} 
\usepackage{ifthen} 
\usepackage{amsfonts}
\usepackage[usenames]{color}
\usepackage{fullpage}
\theoremstyle{plain} 
\newtheorem{thm}{Theorem}

\theoremstyle{definition}

\theoremstyle{remark}
\newtheorem{remark}{Remark}
\newtheorem{ex}{Example}
\newtheorem*{condLP}{Condition LP}
\newtheorem*{condGP}{Condition GP}
\newtheorem*{condS}{Condition S}

\newcommand{\X}{\mathscr{X}}

\newcommand{\E}{\mathsf{E}}

\newcommand{\prob}{\mathsf{P}}

\newcommand{\eps}{\varepsilon}
\renewcommand{\phi}{\varphi}

\renewcommand{\L}{\mathcal{L}}

\newcommand{\nm}{\mathsf{N}}

\newcommand{\unif}{\mathsf{Unif}}
\newcommand{\dir}{\mathsf{Dir}}
\newcommand{\pareto}{\mathsf{Par}}

\title{Empirical priors and posterior concentration rates for a monotone density}
\author{
Ryan Martin\footnote{Department of Statistics, North Carolina State University, email: {\tt rgmarti3@ncsu.edu}}
}
\date{\today}

\begin{document}

\maketitle 

%\doublspacing

\begin{abstract}
%Inference on a monotone density function is an important and challenging problem.  
In a Bayesian context, prior specification for inference on monotone densities is conceptually straightforward, but proving posterior convergence theorems is complicated by the fact that desirable prior concentration properties often are not satisfied.  In this paper, I first develop a new prior designed specifically to satisfy an empirical version of the prior concentration property, and then I give sufficient conditions on the prior inputs such that the corresponding empirical Bayes posterior concentrates around the true monotone density at nearly the optimal minimax rate.  Numerical illustrations also reveal the practical benefits of the proposed empirical Bayes approach compared to Dirichlet process mixtures.

\smallskip

\end{abstract}

\section{Introduction}
\label{S:intro}

Let $X_1,\ldots,X_n$ be iid samples from a density function $f^\star$, supported on the positive half-line, assumed to be monotone non-increasing.  Nonparametric inference on a monotone density has received considerable attention in the literature, dating back to \citet{grenander}, with a wide range of applications \citep[e.g.,][]{robertson.wright.dykstra.1988, groeneboom.jongbloed.book}.  Theoretical properties of estimators have been studied in \citet{prakasarao}, \citet{groeneboom}, and \citet{bala.wellner.2007}, among others, with the behavior of the Grenander estimator at the origin being a now-classical example of inconsistency of the maximum likelihood estimator \citep{woodroofesun} and failure of bootstrap \citep[e.g.,][]{kosorok2008, sen.banerjee.woodroofe.2010}.  

From a Bayesian point of view, constructing a prior and corresponding posterior distribution for the monotone density is at least conceptually straightforward thanks to the mixture representation of \citet{williamson1956}; see Section~\ref{S:prior}.  This makes it possible to construct priors for monotone densities using the standard tools, such as finite mixture models, Dirichlet processes, etc.  However, theoretical analysis of the corresponding posterior distribution is complicated by the fact that, unless the support of $f^\star$ is {\em known}, the usual Kullback--Leibler property \citep[e.g.,][]{ggr1999, wu.ghosal.2008} used to prove posterior convergence results may not be satisfied; in fact, it could be the the Kullback--Leibler divergence of $f$ from $f^\star$ could be infinite for all $f$ in a set of prior probability~1.  Therefore, the general theorems in, e.g., \citet{ggv2000} and \citet{walker2007} cannot be applied.  \citet{salomond2014} worked around this to show, among other things, that the Bayesian posterior distribution based on various mixture priors has concentration rate within a logarithmic factor of the minimax optimal rate, $n^{-1/3}$, with respect to Hellinger or $L_1$ distance.  

A recent trend in the Bayesian literature is asymptotic concentration results for empirical Bayes posteriors; see, e.g., \citet{szabo.vaart.zanten.2013}, \citet{rousseau.szabo.2017}, and \citet{rousseau.etal.eb}.  These papers propose to extend the classical techniques and results to handle the case where the prior involves data in some way, e.g., through a plug-in estimator of a hyperparameter.  However, given that the usual support conditions fail in the problem considered here, even with a fixed prior, it seems unlikely that these new techniques would apply to empirical Bayes monotone density estimation.  \citet{martin.walker.deb}, building on \citet{martin.walker.eb} and \citet{martin.mess.walker.eb}, recently proposed a new empirical Bayes approach, one that constructs the empirical prior specifically so that the desirable posterior concentration rate properties are achieved.  In particular, the empirical prior is designed to satisfy the prior support conditions---a variation on the Kullback--Leibler property---so this approach seems ideally suited for cases, like monotone density estimation, where satisfying the prior support condition is problematic.  Here, in Section~\ref{S:prior}, I will construct a simple and intuitively appealing empirical prior and establish, in Sections~\ref{S:rate}--\ref{S:proofs}, that the corresponding empirical Bayes posterior concentration rate is nearly minimax optimal.  Beyond these desirable theoretical properties, in Section~\ref{S:examples}, I will show that the proposed empirical Bayes approach has a number of practical benefits compared to the Dirichlet process mixture model, including improved computational efficiency and finite-sample performance.  %Proofs are presented in Section~\ref{S:proofs} and some concluding remarks are given in Section~\ref{S:discuss}.  

%The remainder of this paper is organized as follows.  Section~\ref{S:prior} describes the specific empirical prior, suitably centered on an estimator of $f^\star$, and Section~\ref{S:rate} states the main results, namely, conditions on the empirical prior inputs such that the posterior concentrates at the near-optimal minimax rate.  Proofs of the two theorems are presented in Section~\ref{S:proofs} and numerical results are provided in Section~\ref{S:examples}.  The paper concludes with a brief discussion in Section~\ref{S:discuss}. 

\section{An empirical prior}
\label{S:prior}

The starting point here is the representation in \citet{williamson1956} of a monotone density as a scale mixture of uniforms, i.e., for any monotone density density $f$, there exists a mixing distribution $\theta$, supported on a subset of $[0,\infty)$, such that $f = f_\theta$, where 
\[ f_\theta(x) = \int_0^\infty k(x \mid \mu) \, \theta(d\mu), \]
and the kernel $k(x \mid \mu) = \mu^{-1} 1(x \leq \mu)$ is the $\unif(0,\mu)$ density.  From here, a prior for $f$ can be defined by introducing a prior for $\theta$ and using the mapping $\theta \mapsto f_\theta$.  

As is typical, I will model $\theta$ as a (finite) discrete distribution, i.e., 
\[ \theta(d\mu) = \sum_{s=1}^S \omega_s \delta_{\mu_s}(d\mu). \]
This makes $f_\theta$ a finite mixture of uniforms.  For the moment, fix the number of support points $S$.  Then the mixing distribution can be expressed as a finite-dimensional parameter $\theta = (\omega, \mu)$, where $\omega=(\omega_1,\ldots,\omega_S)$ is the vector of mixture weights and $\mu=(\mu_1,\ldots,\mu_S)$ is the corresponding vector of mixture locations.  The theory will require that $S=S_n$ be increasing with $n$ at a suitable rate; see Section~\ref{S:rate}.  For the remainder of this section, I will focus on specifying a prior for $\theta = (\omega, \mu)$, given $S$.  The prior here will be empirical in the sense that it depends on data in a particular way.  

%Regarding the support, there are two cases one can encounter: either the support of the true density $f^\star$ is bounded---with known or unknown upper bound---or the support is unbounded.  In either case, I will consider a prior for the vector $\mu$ that is (basically) restricted to $[t,T]^S$.  If an upper bound on the support is known, then $T$ will be that bound.  If, on the other hand, no upper bound is known, then $T=T_n$ will be suitably increasing in $n$; see Section~\ref{S:rate}.  

The general construction of an empirical prior in \citet{martin.walker.deb} selects an appropriate data-driven center, e.g., the prior mode.  Their motivation is to replace the usual Kullback--Leibler/prior concentration property with an ``empirical'' version.  Towards this, write the likelihood function for $\theta=(\omega, \mu)$, with $S$ fixed, as
\begin{equation}
\label{eq:mixture.lik}
L_n(\theta) = \prod_{i=1}^n \sum_{s=1}^S \omega_s k(X_i \mid \mu_s)
\end{equation}
If $\eps_n$ is the target convergence rate and $\hat\theta$ is a maximizer of the likelihood $L_n$ over a suitable set $\Theta_n$, a {\em sieve}, then \citet{martin.walker.deb} defined
\begin{equation}
\label{eq:Ln}
\L_n = \{\theta \in \Theta_n: L_n(\theta) \geq e^{-dn\eps_n^2} L_n(\hat\theta)\}, \quad d > 0, 
\end{equation}
which is effectively a ``neighborhood'' of $\hat\theta$ in $\Theta_n$, an empirical or data-dependent version of the Kullback--Leibler neighborhood in classical Bayesian nonparametric studies \citep[e.g.,][]{schwartz1965}.  Like in the familiar Bayesian settings, the goal is for the prior to charge $\L_n$ with a sufficient amount of mass; see Condition~LP in Section~\ref{S:proofs}.  But the fact that $\L_n$ is data-dependent means that the prior must also be so, thus, an {\em empirical prior}.  %More specifically, the prior will be centered around $\hat\theta$.  

Here I take the prior mode equal to $\hat\theta = (\hat\omega, \hat\mu)$, a sieve maximum likelihood estimator. The particular sieve is of the form 
\begin{equation}
\label{eq:sieve}
\Theta_n = \{\theta = (\omega, \mu) \in \Delta(S_n) \times [t, T]^{S_n}\}, \quad 0 < t < T < \infty, 
\end{equation}
where $t=t_n$ and $T=T_n$ might also depend on $n$; see Section~\ref{S:rate}.   Then the empirical prior for $\theta=(\omega,\mu)$ I propose here is as follows:
\begin{itemize}
\item $\omega$ and $\mu$ are independent;
\vspace{-2mm}
\item $\omega$ has an $S$-dimensional Dirichlet distribution, $\dir_S(\hat\alpha)$, on the simplex $\Delta(S)$, where $\hat\alpha_s = 1 + c \hat\omega_s$, $s=1,\ldots,S$, and $c=c_n$ is non-stochastic; 
\vspace{-2mm}
\item $\mu_1,\ldots,\mu_S$ are independent with $\mu_s \sim \pareto(\hat\mu_s, \delta)$, a Pareto distribution with scale parameter/lower bound $\hat\mu_s$ and non-stochastic shape parameter $\delta=\delta_n$.  
%$\mu_s \sim \unif(\hat\mu_s, \hat\mu_s + \delta)$, where $\delta=\delta_n$ is a constant to be determined.  
\end{itemize}
It is easy to see that $\hat\theta = (\hat\omega, \hat\mu)$ is the mode. 
%, since $\hat\omega$ is the unique mode for the Dirichlet component.  
The Pareto prior for the components of $\mu$ is convenient because it is conjugate to the uniform mixture kernel.  It is also important for the proofs in Section~\ref{S:proofs} that the prior for $\mu_s$ be supported on $[\hat\mu_s,\infty)$, which is easily arranged with a Pareto distribution.  The to-be-determined constants $(c,\delta)$ control the spread of the prior for $(\omega,\mu)$ around its mode $(\hat\omega,\hat\mu)$.  

%The uniform prior for $\mu$ is only for simplicity and not essential.  The key is that the prior support for $\mu_s$ is bounded, $\hat\mu_s$ is the left end-point, and the density is non-increasing, so other priors would also work.  Since $\delta=\delta_n$ will be vanishing with $n$, the shape of the prior density for $\mu_s$ is largely irrelevant.   

%{\color{red} Finally, the theory in Section~\ref{S:rate} requires some mild restrictions on the optimization that defines $\hat\mu$.  In particular, $\hat\mu_s$ will be the maximum of the likelihood function when restricted to the interval $[X_{(1)}, T]$, where $X_{(1)}$ is the sample minimum.  Given that $T$ is treated as an upper bound on the support of the true density, it is quite reasonable to restrict estimators of the mixture locations to be no more than that bound.  Moreover, it would not be reasonable to have mixture locations smaller than the sample minimum; indeed, {\color{blue} Seregin/Groenboom} has shown that the nonparametric maximum likelihood estimator of $f^\star$ is a finite mixture of the form considered where, with mixture locations no less than $X_{(1)}$, so this is a mild restriction. } 

To summarize, $f$ is modeled as $f_\theta$ and an empirical prior on $f$ is induced by specifying an empirical prior for $\theta$ and using the mapping $\theta \mapsto f_\theta$.  In what follows, $\Pi_n$ will denote the empirical prior for $\theta$ on the sieve \eqref{eq:sieve} as described above.  With a slight abuse of notation, I will also use $\Pi_n$ to denote the corresponding empirical prior for $f=f_\theta$; the meaning should be clear from the context.  Given this prior, the corresponding posterior distribution $\Pi^n$ for $\theta$ is defined as
\begin{equation}
\label{eq:post}
\Pi^n(d\theta) \propto L_n(\theta) \, \Pi_n(d\theta),
\end{equation}
where $L_n(\theta)$ is the likelihood function in \eqref{eq:mixture.lik}.  Again, with a slight abuse of notation, I will also write $\Pi^n$ for the empirical Bayes posterior for the monotone density $f$.

\begin{remark}
\label{re:spread}
The prior support condition eluded to above could be immediately achieved by taking the prior to be degenerate at the $\hat\theta$ that corresponds to Grenander's estimator $\hat f=f_{\hat\theta}$, the nonparametric maximum likelihood estimator.  Of course, the posterior based on this trivial empirical prior is also degenerate at $\hat f$ and, therefore, inherits the concentration rate of Grenander's estimator.  However, achieving the target rate is only a first objective.  By using a non-degenerate prior, the posterior will have spread, leaving open the possibility for uncertainty quantification; see Sections~\ref{S:examples}--\ref{S:discuss}.  
\end{remark}

\begin{remark}
\label{re:clinical}
There is a {\em clinical version} of the model that is perhaps more natural for applications.  In particular, when $n$ is large, the sieve ought to contain the $\theta$ corresponding to Grenander's estimator, so practical applications could dispense with the sieves altogether---which eliminates the need to specify $S$ and $[t,T]$, and to maximize the likelihood over the sieve---and center the prior directly on Grenander's estimator; see Section~\ref{S:examples}.  However, establishing the concentration rate for this clinical version requires control on the mixture support size in Grenander's estimator but, to my knowledge, no such results are available in the literature.  Just like in \citet[][Sec.~4]{ghosalvaart2001}, a reasonable conjecture is that the sieve estimator above is the same as Grenander's, in which case, the clinical version is also covered by Theorems~\ref{thm:bounded}--\ref{thm:unbounded}.  
\end{remark}

\section{Posterior concentration rate}
\label{S:rate}

The previous section described an empirical prior that, when combined with the likelihood via Bayes's formula, leads to a posterior distribution $\Pi^n$ in \eqref{eq:post} that can be used for inference on the monotone density $f$.  But why is this a reasonable approach?  To answer this question, I will provide conditions prior inputs---$c$, $\delta$, $S$, $t$, and $T$---such that the empirical Bayes posterior distribution $\Pi^n$ for $f$ concentrates around the true $f^\star$ at nearly the optimal minimax rate.  Proofs of the two theorems are given in Section~\ref{S:proofs} and finite-sample performance of the posterior is investigated in Section~\ref{S:examples}. %{\color{red} As in Remark~\ref{re:spread}, the goal here is to characterize the ``maximal'' prior spread such that the posterior concentrates at the desired rate.} 

Let $d$ denote the Hellinger or $L_1$ distance.  Then the optimal rate with respect to $d$ is $n^{-1/3}$; see \citet[][Theorem~2.7.5]{vaartwellner1996} and \citet[][Example~3.2]{ggv2000}.  Under certain conditions, this rate can be achieved, within a logarithmic factor, by the nonparametric maximum likelihood estimator \citep[e.g.,][]{birge1989, bala.wellner.2007} and by certain nonparametric Bayesian methods \citep{salomond2014}.  The following theorem establishes the near-optimal concentration rate for the proposed empirical Bayes approach, in the case where $f^\star$ has a bounded support.  

\begin{thm}
\label{thm:bounded}
Let the true density $f^\star$ be monotone non-increasing with support $[0,T^\star]$, with $f^\star(0) < \infty$, and let $\eps_n = (\log n)^{1/3} n^{-1/3}$ be the target rate.  If the prior inputs $(c, \delta, S, t, T) = (c_n, \delta_n, S_n, t_n, T_n)$ satisfy 
\begin{equation}
\label{eq:tuning}
S_n \propto \eps_n^{-1} = n\eps_n^2 (\log n)^{-1}, \quad c_n \propto n \eps_n^{-2}, \quad \text{and} \quad \delta_n \log(T_n / t_n) \lesssim \log n, 
%\delta_n^{-1} T_n \propto n, 
\end{equation} 
%\begin{align*}
%c & = c_n = n \eps_n^{-2} \\ %= n^{4/3} (\log n)^{-1/3} \\
%\delta & = \delta_n \propto n^{-1} \\
%S & =S_n \propto n \eps_n^2 (\log n)^{-1}, %= n^{1/3} (\log n)^{-1/3},
%\end{align*}
where $T_n$ and $\delta_n$ are non-decreasing and $t_n$ is non-increasing, then there exists a constant $M > 0$ such that the posterior distribution $\Pi^n$ satisfies 
\begin{equation}
\label{eq:rate}
\E_{f^\star}\bigl[ \Pi^n(\{f: d(f^\star, f) > M \eps_n\}) \bigr] \to 0, \quad n \to \infty. 
\end{equation}
%In particular, if $T^\star$ is known, then one can take $T_n \equiv T^\star$ and $\delta_n \propto n^{-1}$.  
\end{thm}

%\begin{proof}
%See Section~\ref{S:proofs}.
%\end{proof}

The conclusion here is similar to that in Theorem~1 of \citet{salomond2014}.  Indeed, the concentration rate above is the same as that obtained by a suitable Dirichlet process mixture model, which is minimax optimal up to the logarithmic factor.  %However, as will be made clear in Section~\ref{S:examples}, in addition to the optimal theoretical concentration rate, there are some practical benefits to the proposed empirical Bayes approach compared to the fully Bayesian Dirichlet process mixture model formulation.  

The latter condition in \eqref{eq:tuning} deserves some explanation.  Since the upper bound $T^\star$ is finite, any sufficiently large $T_n$ would suffice for estimating $f^\star(x)$ for $x$ near $T^\star$.  Similarly, since a draw $f$ from the proposed prior satisfies $f(0) \leq t_n^{-1}$, any $t_n$ less than $f^\star(0)^{-1} > 0$ would suffice for estimating $f^\star(x)$ for $x$ near 0.  Therefore, $\log(T_n/t_n)$ could be very slowly increasing, or even bounded, which means $\delta_n$ can grow as fast as order $\log n$.  Intuitively, slowly increasing $\log(T_n/t_n)$ indicates some certainty about the support of $f^\star$, in which case, a larger $\delta_n$ and, hence, a smaller Pareto variance, is reasonable.  On the other hand, if the support is uncertain, one could take $T_n/t_n$ to be polynomial in $n$, in which case $\delta_n$ must be bounded and, hence, the Pareto variance must be bounded away from zero.

%However, the setup here is arguably simpler and the empirical Bayes computations---via Gibbs sampling---are faster and easier than for the corresponding Dirichlet process mixture model; see Section~\ref{S:examples}; in a real application, where $n$ is fixed, approximation of the above empirical Bayes posterior is easy via Gibbs sampling.  Salomond's setup, on the other hand, assumes either a Dirichlet process mixture prior for $f$ or a finite mixture with a prior on the number of components.  These latter models are common, but apparently nothing is gained theoretically from their additional computational complexity compared to this simple empirical Bayes approach.  

The only serious assumption on $f^\star$ in Theorem~\ref{thm:bounded} is that the support is bounded.  It turns out that this bounded-support condition can be replaced by a condition on the tails of $f^\star$.  Condition~C4 in \citet{salomond2014} states:
\begin{equation}
\label{eq:tail}
\text{there exists $b, r > 0$ such that $f^\star(x) \leq e^{-b x^r}$ for all large $x$}. 
\end{equation}
The next result is analogous to Theorem~2 in \citet{salomond2014}.  

\begin{thm}
\label{thm:unbounded}
Let the true density $f^\star$ be monotone non-increasing, with $f^\star(0) < \infty$, whose support is $[0, \infty)$.  Assume that $f^\star$ satisfies \eqref{eq:tail} for a given $r$, and set the target rate equal to $\eps_n = (\log n)^{1/3 + 1/r} n^{-1/3}$.  Let the prior inputs be as in \eqref{eq:tuning}, but with $T_n \gtrsim (\log n)^{1/r}$.  Then the conclusion \eqref{eq:rate} of Theorem~\ref{thm:bounded} holds with the modified rate $\eps_n$.  
\end{thm}

%\begin{proof}
%See Section~\ref{S:proofs}
%\end{proof}

Note that the rate in the unbounded support case is slightly slower, by only a logarithmic factor, than in the bounded support case.  Actually, Theorem~\ref{thm:bounded} can be viewed as a special case of this Theorem~\ref{thm:unbounded}: if the support is bounded, then \eqref{eq:tail} holds for ``$r=\infty$'' and so the rate in Theorem~\ref{thm:unbounded} agrees with that in Theorem~\ref{thm:bounded}.  The only subtlety is that the choice of $T=T_n$ depends on $r$, a feature of $f^\star$, which is typically unknown in real examples.  If one is willing to assume a positive lower bound $r_0$ on $r$, then $T=T_n$ can be chosen with $r=r_0$.  It is not known if the rate in Theorem~\ref{thm:unbounded} is optimal so, even though working with a lower bound on $r_0$---corresponding to a larger $T_n$---will slow down the rate slightly, there is no practical difference compared to the rate with the true $r$.  \citet{salomond2014} does not say whether his Theorem~2 requires knowledge of the tail exponent $r$, but the tails of the Dirichlet process base measure generally affect the posterior concentration rates \citep[e.g.,][Theorem~5.1]{ghosalvaart2001}.  So if a target rate depending on $r$ is to be achieved, then this requires some $r$-dependent condition on the base measure and, therefore, to check this condition, $r$ must be known.

\section{Proofs}
\label{S:proofs}

\subsection{General strategy}
\label{SS:strategy}

\citet{martin.walker.deb} proposed a general strategy for constructing empirical priors such that the corresponding posterior distribution has the desired concentration properties.  Their Theorem~1 lists three general conditions that include assumptions about the prior concentration, one local and one global, as well as an assumption about the approximation properties of the sieve.  I will summarize these conditions in the context of iid data as being considered here.  Let $\eps_n$ be the target rate.  

\begin{condS}
There exists a $\theta^\dagger = \theta_n^\dagger$ in the sieve $\Theta_n$ such that 
\[ \max\Bigl\{ \int \Bigl( \log \frac{f^\star}{f_{\theta^\dagger}} \Bigr) f^\star \,dx, \, \int \Bigl( \log \frac{f^\star}{f_{\theta^\dagger}} \Bigr)^2 f^\star \, dx \Bigr\} \leq \eps_n^2. \]
%\[ \max\bigl\{ K(f^\star, f_{\theta^\dagger}), V(f^\star, f_{\theta^\dagger}) \bigr\} \leq \eps_n^2, \]
%where $K$ and $V$ are the Kullback--Leibler and corresponding second moment
\end{condS}

\begin{condLP}
For a given $d > 0$ define $\L_n$ as in \eqref{eq:Ln}.  Then there exists a constant $C > 0$ such that the empirical prior $\Pi_n$ satisfies 
\[ \liminf_{n \to \infty} e^{C n \eps_n^2} \Pi_n(\L_n) > 0, \quad \text{with $\prob_{f^\star}$-probability 1}. \]
\end{condLP}

\begin{condGP}
Let $\pi_n$ be the density function for $\theta$ under the empirical prior.  For a constant $p > 1$, there exists $K > 0$ such that 
\begin{equation}
\label{eq:gp}
\int_{\Theta_n} \bigl[ \E_{f^\star}\{ \pi_n(\theta)^p \} \bigr]^{1/p} \,d\theta \leq e^{K n \eps_n^2}. 
\end{equation}
\end{condGP}

\subsection{Proof of Theorem~\ref{thm:bounded}}

I will begin by checking Condition~LP.  For fixed $S$, if $\theta=(\omega, \mu)$ denotes the mixture weights and locations, respectively, then the likelihood function $L_n(\theta)$ in \eqref{eq:mixture.lik} for the discrete mixture model can be expressed as
\begin{align*}
L_n(\theta) & = \sum_{(n_1,\ldots,n_S)} \omega_1^{n_1} \cdots \omega_S^{n_S} \sum_{(s_1,\ldots,s_n)} \prod_{s=1}^S \prod_{i: s_i = s} k(X_i \mid \mu_s) \\
& = \sum_{(n_1,\ldots,n_S)} \omega_1^{n_1} \cdots \omega_S^{n_S} \sum_{(s_1,\ldots,s_n)} \prod_{s=1}^S \mu_s^{-n_s} 1(\mu_s \geq \hat X_s), 
\end{align*}
where $\hat X_s = \max_{i: s_i=s} X_i$ is the largest of the $n_s$ $X$ values in category $s$, relative to the partition determined by $(s_1,\ldots,s_n)$ with frequency table $(n_1,\ldots,n_S)$; the inner- and outer-most sums are over all such partitions and frequency tables, respectively.  Since the prior for $\mu_s$ is supported on $[\hat\mu_s,\infty)$, I can lower-bound the likelihood by 
\[ L_n(\theta) \geq \sum_{(n_1,\ldots,n_S)} \omega_1^{n_1} \cdots \omega_S^{n_S} \sum_{(s_1,\ldots,s_n)} \prod_{s=1}^S \mu_s^{-n_s} 1(\hat\mu_s \geq \hat X_s). \]
The prior has $\omega$ and $\mu$ independent, and $\mu_1,\ldots,\mu_S$ independent, so 
\[ \E\{L_n(\theta)\} \geq \sum_{(n_1,\ldots,n_S)} \E(\omega_1^{n_1} \cdots \omega_S^{n_S}) \sum_{(s_1,\ldots,s_n)} \prod_{s=1}^S \E(\mu_s^{-n_s}) 1(\hat\mu_s \geq \hat X_s), \]
where expectation is with respect to the prior for $\theta=(\omega, \mu)$.  The proof of Proposition~2 in \citet{martin.walker.deb} gives a bound for the first expectation, i.e., 
\[ \E(\omega_1^{n_1} \cdots \omega_S^{n_S}) \geq \frac{\Gamma(c + S)c^n}{\Gamma(c + S + n)} \hat\omega_1^{n_1} \cdots \hat\omega_S^{n_S}. \]
For the Pareto prior, $\pareto(\hat\mu_s, \delta)$, on $\mu_s$, we have 
%if $n_s > 1$, then we have 
\[ \E(\mu_s^{-n_s}) = \int_{\hat\mu_s}^\infty \frac{\delta \hat\mu_s^\delta}{\mu_s^{\delta + n_s +1}} \,d\mu_s = \frac{\delta}{\delta + n_s} \hat\mu_s^{-n_s} \geq \frac{1}{1 + n\delta^{-1}} \hat\mu_s^{-n_s}. \]
%\[ \E(\mu_s^{-n_s}) = \frac{1}{\delta} \int_{\hat\mu_s}^{\hat\mu_s + \delta} \mu_s^{-n_s} \,d\mu_s = \frac{1}{n_s - 1} \frac{1}{\delta} \frac{1}{\hat\mu_s^{n_s - 1}} \Bigl\{ 1 - \Bigl( \frac{\hat\mu_s}{\hat\mu_s + \delta} \Bigr)^{n_s - 1} \Bigr\}. \]
%But $\delta$ is small, i.e., $\delta = \delta_n \to 0$, so the ratio of the difference $\{\cdots\}$ to $\delta$ is approximately the derivative of $\delta \mapsto -(\frac{\hat\mu_s}{\hat\mu_s + \delta} )^{n_s-1}$ at $\delta = 0$, which is $(n_s-1) / \hat\mu_s^{n_s}$.  Therefore, 
%\[ \E(\mu_s^{-n_s}) \geq \tfrac12 \hat\mu_s^{-n_s}, \quad \text{all large $n$}. \]
%The fraction $\frac12$ is not essential, it can be replaced by any other number in $(0,1)$.  The same bound holds with $n_s=1$, but the calculus is a bit different/easier.  
Since $\delta=\delta_n$ is non-vanishing, the first term in the lower bound is at least $O(n^{-1})$.  Plugging this bound back into the above expectation gives
\[ \E\{L_n(\theta)\} \geq \frac{\Gamma(c + S)c^n}{\Gamma(c + S + n)} e^{-S \log n} L_n(\hat\theta), \quad \text{all large $n$}. \]
As in Proposition~2 of \citet{martin.walker.deb}, if $S=S_n$ is of the order $n\eps_n^2 (\log n)^{-1}$ and $c = n \eps_n^{-2}$, then there exists a constant $K > 0$ such that 
\[ \frac{\Gamma(c + S)c^n}{\Gamma(c + S + n)} \geq e^{-K n \eps_n^2}. \]
Similarly, for the second term, since $S \log n$ is of the order $n\eps_n^2$, we can conclude that, for a suitable constant $D > 0$, 
\[ \E\{L_n(\theta)\} / L_n(\hat\theta) > e^{-D n \eps_n^2} \]
which, according to the argument in the proof of Proposition~3 in \citet{martin.walker.deb}, implies Condition~LP.  

Next, I check Condition~GP.  As a first step, we have that the density for the Dirichlet prior on $\omega$ is uniformly upper bounded by $(c+S)^{c+S+1/2}c^{-c}$, which does not depend on data.  For the prior on $\mu$, the density function is upper bounded by 
\[ \delta T^\delta \mu_s^{-(\delta + 1)} \, 1(\mu_s \geq t), \]
%\[ \delta^{-S} \prod_{s=1}^S 1(\mu_s \leq T + \delta), \]
which is also free of data.  Then, for any $p > 1$, the integral \eqref{eq:gp} is bounded by 
\begin{equation}
\label{eq:bound1}
\frac{(c+S)^{c+S+1/2}}{c^c} %\Bigl( 1 + \frac{T}{\delta} \Bigr)^S 
\, \Bigl(\frac{T}{t} \Bigr)^{\delta S}. 
\end{equation}
With $S=S_n$ of the order $n\eps_n^2 (\log n)^{-1}$ and 
%$\delta^{-1}T=\delta_n^{-1}T_n$ of the order $n$, 
$\delta \log(T/t)$ of order $\log n$, it follows that the second term in \eqref{eq:bound1} is bounded by $e^{A n \eps_n^2}$.  Similarly, for $c = n \eps_n^2$ and $S=S_n$ of the order $n \eps_n^2 (\log n)^{-1}$, \citet{martin.walker.deb} showed that the first term in \eqref{eq:bound1} is also $e^{B n \eps_n^2}$ so, altogether, the relevant integral is bounded by $e^{C n \eps_n^2}$, hence Condition~GP.  

Finally, note that, for the case of $f^\star$ with bounded support $[0,T^\star]$ and $f^\star(0) < \infty$, Condition~S on the sieve follows from Lemma~11 in \citet{salomond2014}.  My sieve has a lower bound $t > 0$ but, if it is small enough, then it does not affect Salomond's calculations since there is no benefit to having a mixture location smaller than $f^\star(0)^{-1} > 0$.  Having checked the three conditions in Section~\ref{SS:strategy}, the conclusion of Theorem~\ref{thm:bounded} follows from the general results in \citet{martin.walker.deb}.

\subsection{Proof of Theorem~\ref{thm:unbounded}}

For a given $T$, possibly depending on $n$, write $f_T^\star$ for the normalized version of $f^\star$ to the interval $[0,T]$, i.e., $f_T^\star(x) = f^\star(x) / F^\star(T)$, where $F^\star$ is the distribution function corresponding to $f^\star$.  Without loss of generality, let $d$ denote the $L_1$ distance.  Then the triangle inequality implies that, for any density $f$, 
\[ d(f^\star, f) \leq d(f^\star, f_T^\star) + d(f_T^\star, f). \]
Moreover, a simple calculation shows that 
\[ d(f^\star, f_T^\star) = \frac12 \int_0^\infty |f^\star(x) - f_T^\star(x)| \,dx = 1 - F^\star(T) . \]
Under the condition \eqref{eq:tail} on the density $f^\star$, it is easy to check that the tail probability $1-F^\star(T) \lesssim \Gamma(r^{-1}, T^r)$, where $\Gamma(s,t)$ is the upper incomplete gamma function, i.e., $\Gamma(s,t) = \int_t^\infty y^{s-1} e^{-y} \,dy$.  From the well-known asymptotic behavior of this gamma function, it follows that if $T_n \gtrsim (\log n)^{1/r}$, then 
\begin{equation}
\label{eq:cdf.bound}
1-F^\star(T_n) \lesssim n^{-1}\eps_n, 
\end{equation}
where $\eps_n = (\log n)^{1/3 + 1/r} n^{-1/3}$ as in the statement of the theorem; see, also, page 1390 in \citet{salomond2014}.  Therefore, for any $f$, 
\begin{align*}
d(f^\star, f) > M \eps_n & \implies d(f^\star, f_T^\star) + d(f_T^\star, f) > M \eps_n \\
& \implies d(f_T^\star, f) > M \eps_n - d(f^\star, f_T^\star) = M \eps_n - \{1 - F^\star(T_n)\} \\
& \implies d(f_T^\star, f) > (M/2) \, \eps_n, \quad \text{say}.
\end{align*}
This effectively converts the problem into one with bounded support.  To see this, define the two sets of densities 
\[ A_n = \{f: d(f^\star, f) > M \eps_n\} \quad \text{and} \quad B_n = \{f: d(f_{T_n}^\star, f) > (M / 2) \, \eps_n\}. \]
Then the argument above implies that $A_n \subseteq B_n$ which, in turn, implies 
\[ \Pi^n(A_n) \leq \Pi^n(B_n). \]
Define the event $\X_n = \{(x_1,\ldots,x_n) \in [0,\infty)^n: x_{(n)} \leq T_n\}$, where $x_{(n)} = \max_i x_i$.  Based on the bound in \eqref{eq:cdf.bound}, it is easy to check that $\prob_{f^\star}(\X_n^c) = o(1)$.  Next, write 
\[ \Pi^n(B_n) = \Pi^n(B_n) 1(\X_n) + \Pi^n(B_n) 1(\X_n^c) \leq \Pi^n(B_n) 1(\X_n) + 1(\X_n^c), \]
so that 
\[ \E_{f^\star}\{\Pi^n(B_n)\} \leq \E_{f^\star}\{\Pi^n(B_n) 1(\X_n)\} + o(1). \]
The expectation on the right-hand side can be rewritten as 
\[ \E_{f^\star}\{\Pi^n(B_n) 1(\X_n)\} = \E_{f^\star}\{\Pi^n(B_n) \mid \X_n\} \prob_{f^\star}(\X_n). \]
Since $\prob_{f^\star}(\X_n) \to 1$, it remains to deal with the conditional expectation.  The key observation is that the conditional distribution of $(X_1,\ldots,X_n)$, given $\X_n$, is iid $f_{T_n}^\star$ and, therefore, 
\begin{equation}
\label{eq:conditioning}
\E_{f^\star}\{\Pi^n(B_n) \mid \X_n\} = \E_{f_{T_n}^\star}\{\Pi^n(B_n)\}, 
\end{equation}
hence, the claim that this is effectively a bounded support problem.  Moreover, all of the work in checking Conditions~LP, GP, and S in the proof of Theorem~\ref{thm:bounded} above applies here with $f^\star$ replaced by $f_{T_n}^\star$ and the modified rate.  Since the general results in \citet{martin.walker.deb} do not require that the ``true parameter'' be fixed, we can conclude that the right-hand side of \eqref{eq:conditioning} is $o(1)$ and, therefore, 
\[ \E_{f^\star}\{\Pi^n(A_n)\} \leq \{1 + o(1)\} \E_{f_{T_n}^\star} \{\Pi^n(B_n)\} + o(1) \to 0, \]
as was to be shown.  This completes the proof of Theorem~\ref{thm:unbounded}.

\section{Numerical results}
\label{S:examples}

\subsection{Computation}
\label{SS:comp}

Here I will focus on the {\em clinical version} of the empirical prior suggested in Remark~\ref{re:clinical}.  That is, this version of the prior is centered on Grenander's estimator, which is obtained from the output provided by the function {\tt grenander} in the R package ``fdrtool'' \citep{fdrtool}; more precisely, the clinical empirical prior is centered on the $\hat\theta=(\hat\omega, \hat\mu)$ for which $f_{\hat\theta}$ is Grenander's estimator.  For moderate to large $n$, there ought to be no difference between this and the empirical prior in Section~\ref{S:prior}, but the former has two clear advantages: first, there is no need to specify $S$ or the interval $[t,T]$; second, maximizing the likelihood over the sieve $\Theta_n$ is not so easy, but there is already efficient software available for computing Grenander's estimator.  

Since $S$ and $[t,T]$ are taken care of by the clinical formulation, it only remains to specify $c=c_n$ and $\delta=\delta_n$.  In what follows, I will use 
\begin{equation}
\label{eq:hyper}
c = 0.01 \frac{n^{5/3}}{(\log n)^{2/3}} \quad \text{and} \quad \delta = \frac{\log n}{20}, 
\end{equation}
a choice guided by the conditions in Theorems~\ref{thm:bounded}--\ref{thm:unbounded}.  The simulation experiments below suggest that this choice works well in various cases, though more work would be needed to determine if these values are ``good'' in any general sense; see Section~\ref{S:discuss}.  

I will compare the empirical Bayes results with those from a Bayesian Dirichlet process mixture model.  I take the Dirichlet process precision parameter as 1 and the base measure as an inverse gamma with data-driven choice of shape and scale parameters, $a$ and $b$.  That is, I choose $a$ and $b$ so that the base measure mean and variance match the mean and variance of the data; this is to ensure that the base measure center and spread are consistent with the data.  Other base measures could also be considered, such as gamma.  A Pareto base measure would be a natural choice, given its conjugacy to the uniform kernel, but this does not satisfy the tail conditions in \citet{salomond2014} needed to achieve the nearly-$n^{-1/3}$ concentration rate.  Conjugacy would allow for very convenient posterior sampling via the slice sampler in \citet{walker2007.slice} and \citet{kalli.griffin.walker.2011}.  But for other base measures, such as gamma and inverse gamma, rejection sampling is required in every MCMC iteration, so the computation is slower than with a Pareto base measure.  Note that the proposed empirical Bayes approach {\em can} employ the computationally efficient Pareto prior since the strategic data-driven centering makes the tails irrelevant to the concentration rate properties.  
%Pareto distribution with scale/lower bound parameter as $\min_i X_i$ and the shape parameter such that the data range, $[\min_i X_i, \max_i X_i]$, is assigned prior probability 0.90; throughout, I set the Dirichlet process precision parameter at 1.  Other base measures could be considered, but the Pareto's conjugate form to the uniform kernel is convenient because it handles the varying support constraints much more easily and efficiently than, say, an inverse gamma base measure.  
%The particular algorithm I employed is the slice sampler in \citet{walker2007.slice} and \citet{kalli.griffin.walker.2011}.  
R codes for both the Bayes and empirical Bayes methods are available at \url{https://www4.stat.ncsu.edu/~rmartin/}.

\subsection{Simulated data}

\begin{figure}[t]
\begin{center}
\subfigure[Empirical Bayes]{\scalebox{0.6}{\includegraphics{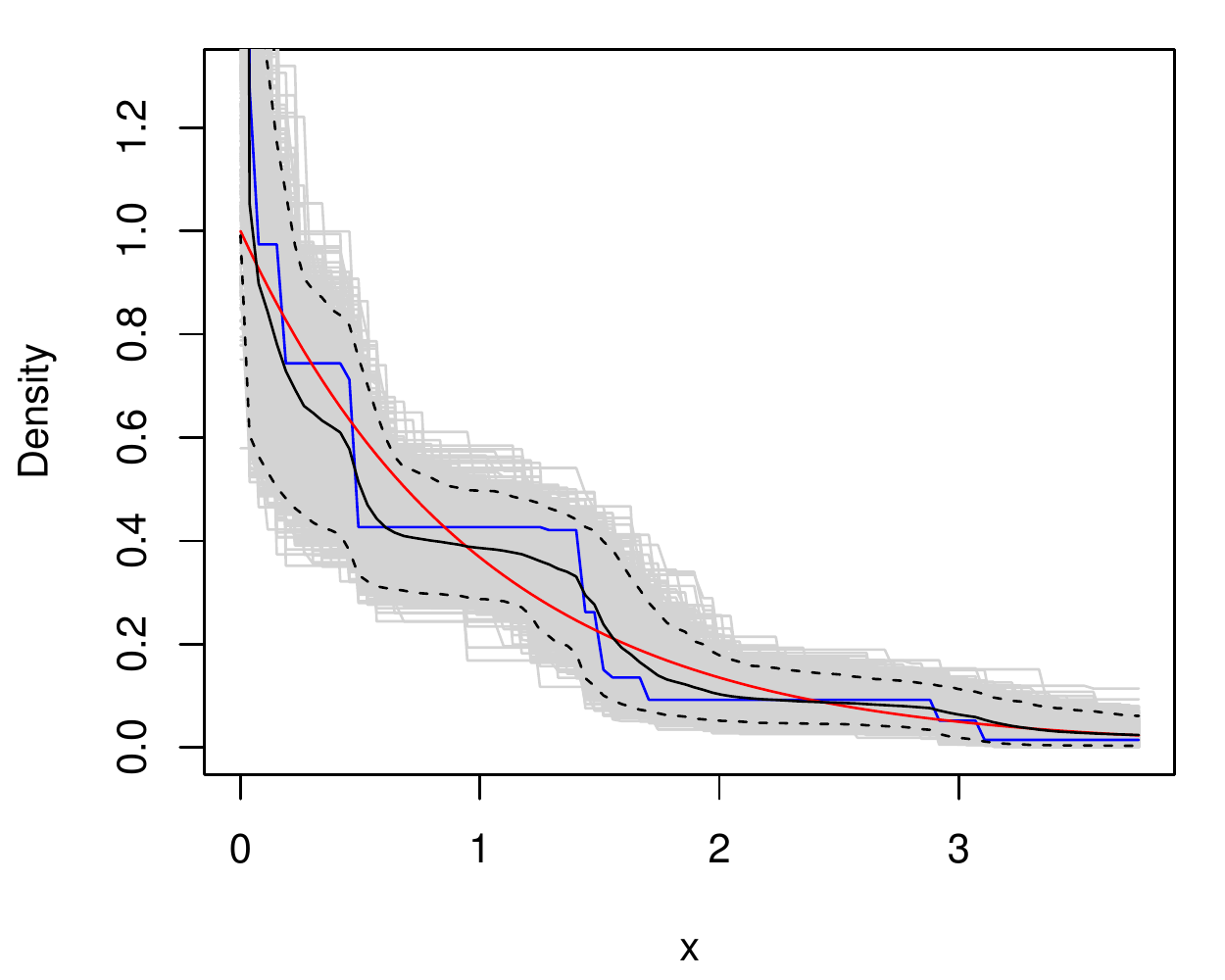}}}
\subfigure[Dirichlet process mixture]{\scalebox{0.6}{\includegraphics{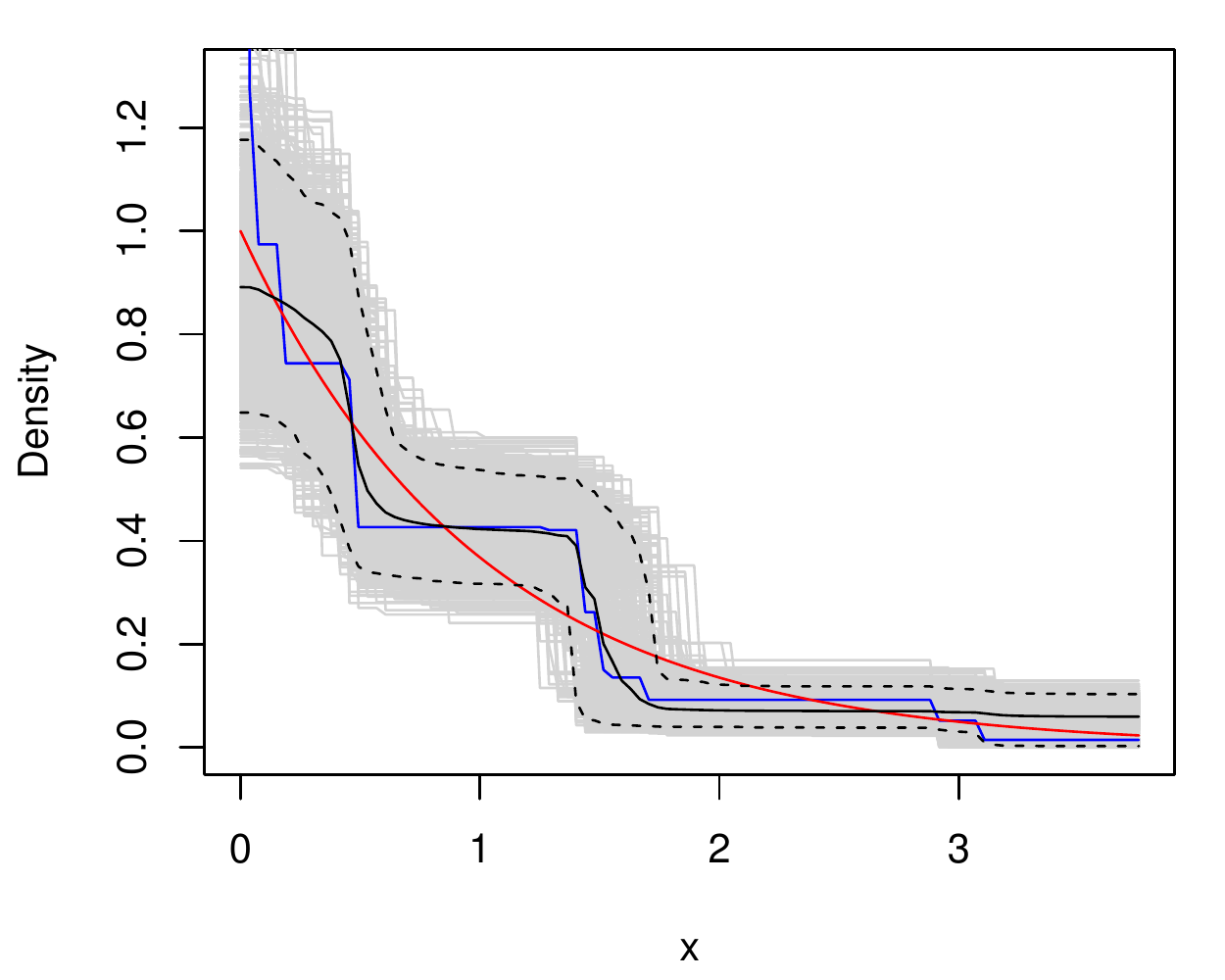}}}
\end{center}
\caption{Exponential example.  Plots of posterior samples (gray), the posterior mean (black), 95\% credible band (dashed), the Grenander estimator (blue), and $f^\star$ (red).}
\label{fig:simulation}
\end{figure}

Consider the case where $f^\star(x) = e^{-x}$, $x \geq 0$, is an exponential density.  Figure~\ref{fig:simulation} plots posterior samples from the respective posterior distributions, the posterior mean density, the 95\% posterior credible band, the Grenander estimator, and the truth.  Both posteriors generally follow the Grenander estimator, and both have roughly the same spread but, at least in this case, the empirical Bayes posterior is less rough and better captures the shape of the truth.  The credible bands for both posteriors cover the truth at most $x$ values, but the Dirichlet process mixture misses in a few places because of its roughness.  Another difference, not apparent from the plots, is that the empirical Bayes computations take about 50\% less time than the Dirichlet process mixture. 

Next I present a more detailed investigation into the concentration properties of the two posteriors.  Specifically, I will consider the coverage properties of the 95\% credible bands at select values of $x$.  That is, for a fixed $x$, I get posterior samples of $f(x)$ and return the 0.025 and 0.975 quantiles as the 95\% posterior credible interval.  I repeat this process for 500 data sets and report the (Monte Carlo approximation of the) coverage probability and expected lengths for the two posteriors.  Table~\ref{tab:expsim} reports (Monte Carlo approximations of) the coverage probability and mean lengths for the 95\% credible intervals based on 500 samples from an exponential distribution, where $f^\star(x) = e^{-x}$, $x \geq 0$.  Table~\ref{tab:hnmsim} reports the same but for a half-normal distribution, where $f^\star(x) = 2 \, \nm(x \mid 0, 1)$, $x \geq 0$. 

The key observation across all the different setting is that the 95\% posterior credible intervals from the Dirichlet process mixture model tend to under-cover, i.e., the coverage probability is less than the nominal 0.95, sometimes much less, while those from the empirical Bayes approach tend to over-cover.  However, the higher coverage of the empirical Bayes intervals is not a result of being wider on average; in fact, the empirical Bayes credible intervals tend to be {\em shorter} than the Bayesian competitor's.  Coverage properties of nonparametric methods is a delicate matter and a detailed investigation is beyond the scope of this paper, but these results strongly suggest that the proposed empirical Bayes procedure---despite its ``double-use'' of the data---does not follow the data too closely and, as suggested in Remark~\ref{re:spread}, may provide valid uncertainty quantification.  

\begin{table}[t]
\begin{center}
\begin{tabular}{cccccc}
\hline
& & \multicolumn{2}{c}{Coverage Prob.} & \multicolumn{2}{c}{Mean Length} \\
\cline{3-4} \cline{5-6} 
$n$ & $x$ & EB & DPM & EB & DPM \\
\hline
100 & 0.5 & 0.943 & 0.987 & 0.371 & 0.459 \\
& 1.0 & 0.984 & 0.965 & 0.270 & 0.338 \\
& 2.0 & 0.986 & 0.938 & 0.150 & 0.185 \\
& 3.0 & 0.990 & 0.891 & 0.086 & 0.098 \\
\hline 
200 & 0.5 & 0.934 & 0.981 & 0.285 & 0.371 \\
& 1.0 & 0.962 & 0.955 & 0.219 & 0.285 \\
& 2.0 & 0.983 & 0.914 & 0.120 & 0.154 \\
& 3.0 & 0.970 & 0.842 & 0.068 & 0.082 \\
\hline
\end{tabular}
\end{center}
\caption{Exponential example---coverage probability and mean length of the 95\% posterior credible regions for the empirical Bayes (EB) and Dirichlet process mixture model (DPM), for several values of $x$ and $n$.}
\label{tab:expsim}
\end{table}

\begin{table}[t]
\begin{center}
\begin{tabular}{cccccc}
\hline
& & \multicolumn{2}{c}{Coverage Prob.} & \multicolumn{2}{c}{Mean Length} \\
\cline{3-4} \cline{5-6} 
$n$ & $x$ & EB & DPM & EB & DPM \\
\hline
100 & 0.5 & 0.896 & 0.974 & 0.313 & 0.370 \\
& 1.0 & 0.969 & 0.961 & 0.291 & 0.361 \\
& 2.0 & 0.971 & 0.897 & 0.178 & 0.205 \\
& 3.0 & 0.986 & 0.636 & 0.054 & 0.043 \\
\hline
200 & 0.5 & 0.865 & 0.991 & 0.242 & 0.294 \\
& 1.0 & 0.964 & 0.954 & 0.235 & 0.304 \\
& 2.0 & 0.956 & 0.885 & 0.137 & 0.166 \\
& 3.0 & 0.951 & 0.717 & 0.043 & 0.045 \\
\hline
\end{tabular}
\end{center}
\caption{Half-normal example---same results as in Table~\ref{tab:hnmsim}.}
\label{tab:hnmsim}
\end{table}

\ifthenelse{1=1}{}{
> o.exp.100 <- cvg.compare(1, 200, 100, 2000, 0.01, 20); print(o.exp.100)
       x eb.cvg dp.cvg     eb.len     dp.len eb.len.med dp.len.med
[1,] 0.0  0.955  0.810 9.38103309 0.52528569 2.97308812 0.49697929
[2,] 0.5  0.945  0.985 0.37058961 0.45906381 0.36298185 0.45528183
[3,] 1.0  0.985  0.965 0.27029776 0.33849662 0.26491170 0.32528564
[4,] 2.0  0.985  0.935 0.14954669 0.18502112 0.14608297 0.17504969
[5,] 3.0  0.990  0.890 0.08600687 0.09769028 0.08383119 0.08991125

> o.exp.200 <- cvg.compare(1, 200, 200, 2000, 0.01, 20); print(o.exp.200)
       x eb.cvg dp.cvg      eb.len     dp.len eb.len.med dp.len.med
[1,] 0.0  0.870  0.565 25.10167953 0.37743097 2.74706845 0.35777078
[2,] 0.5  0.935  0.980  0.28495878 0.37091416 0.27244490 0.35947146
[3,] 1.0  0.960  0.955  0.21934588 0.28509377 0.21352195 0.27993690
[4,] 2.0  0.985  0.915  0.12006742 0.15410351 0.11800884 0.14459434
[5,] 3.0  0.970  0.840  0.06753994 0.08162171 0.06619479 0.07565278

> o.hnm.100 <- cvg.compare(2, 200, 100, 2000, 0.01, 20); print(o.hnm.100)
       x eb.cvg dp.cvg     eb.len    dp.len eb.len.med dp.len.med
[1,] 0.0  0.870  0.975 15.5818394 0.4661701 2.28964769 0.44757330
[2,] 0.5  0.895  0.975  0.3130285 0.3695994 0.30649681 0.36055540
[3,] 1.0  0.970  0.960  0.2913823 0.3614509 0.28060416 0.34774879
[4,] 2.0  0.970  0.895  0.1775161 0.2053701 0.17430277 0.19529811
[5,] 3.0  0.985  0.635  0.0535107 0.0430474 0.05304542 0.03828106

> o.hnm.200 <- cvg.compare(2, 200, 200, 2000, 0.01, 20); print(o.hnm.200)
       x eb.cvg dp.cvg      eb.len     dp.len eb.len.med dp.len.med
[1,] 0.0  0.835  0.990 39.13048172 0.32175141 2.33909254 0.31056913
[2,] 0.5  0.865  0.990  0.24174820 0.29383375 0.23581424 0.28729490
[3,] 1.0  0.965  0.955  0.23525560 0.30390177 0.22520438 0.29365285
[4,] 2.0  0.955  0.885  0.13719649 0.16600999 0.13554070 0.16757394
[5,] 3.0  0.950  0.715  0.04280447 0.04472443 0.04325217 0.04645635
}

\subsection{Real data}

%Here I provide two illustrations of the posterior based on the empirical prior described in Section~\ref{S:prior}.  Again, comparison will be made with the results obtained from the Dirichlet process mixture model as described above.  

\begin{figure}[t]
\begin{center}
\subfigure[Data and empirical Bayes]{\scalebox{0.6}{\includegraphics{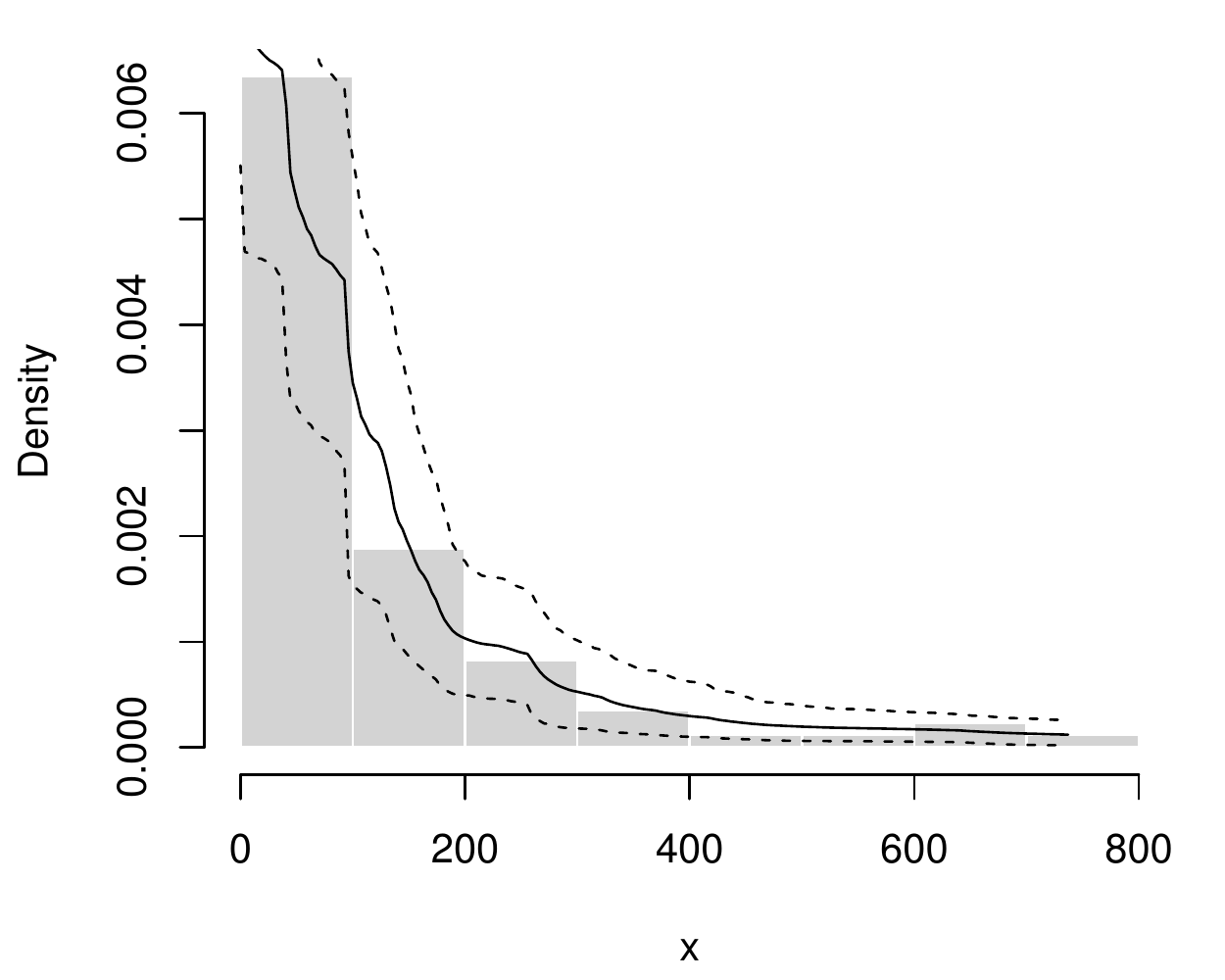}}}
\subfigure[Data and Dirichlet process mixture]{\scalebox{0.6}{\includegraphics{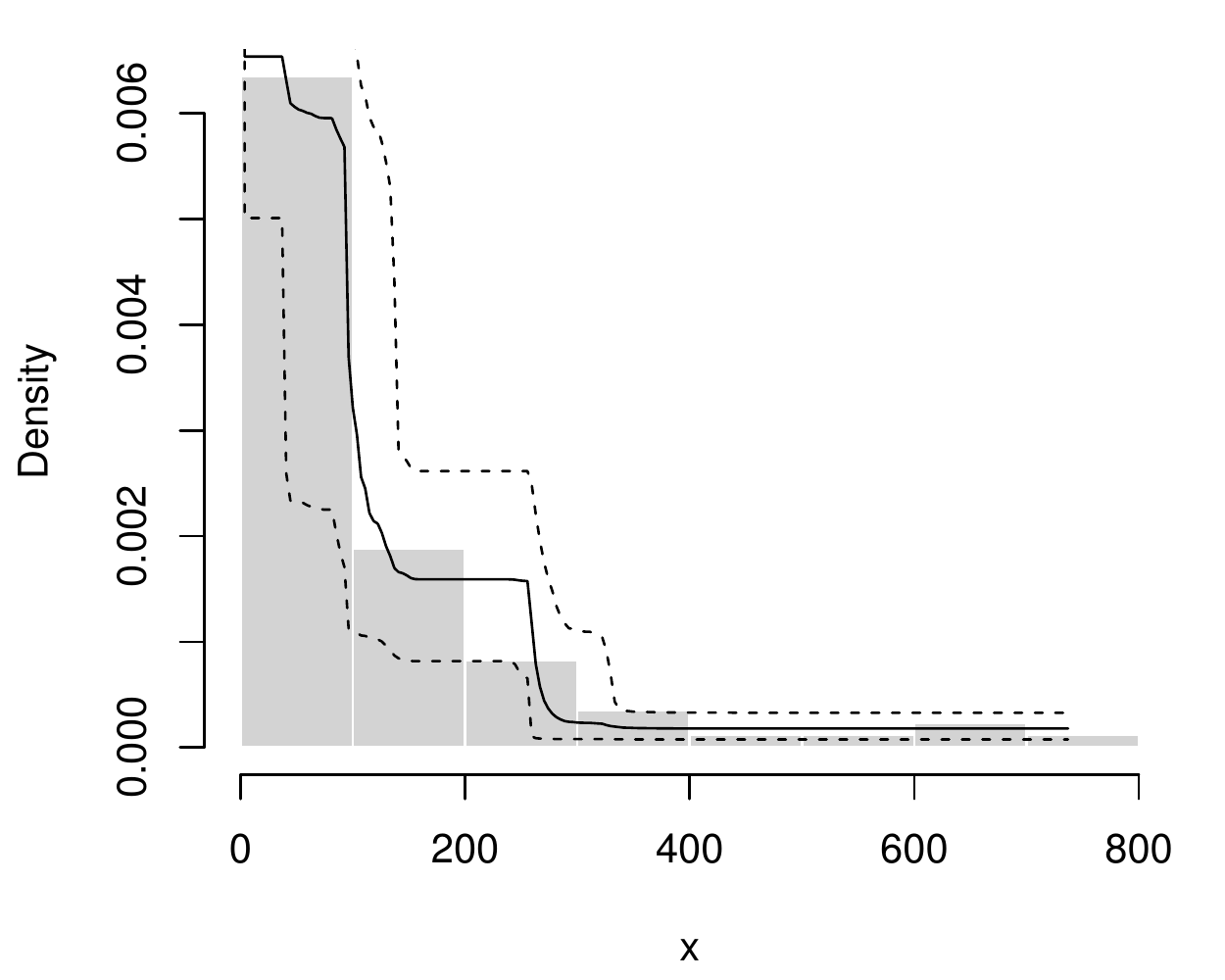}}}
\end{center}
\caption{Example~\ref{ex:silverman}, with $n=86$---data histogram, posterior mean (solid), and 95\% credible band (dashed).}
\label{fig:silverman}
\end{figure}

\begin{ex}
\label{ex:silverman}
\citet[][Table~2.1]{silverman} presents data on the lengths of psychiatric treatment undergone by $n=86$ patients used as controls in a study of suicide risks.  Figure~\ref{fig:silverman} shows the data histogram, the posterior mean density, and the 95\% credible bands for the two estimators, based on the same settings described in Section~\ref{SS:comp}.  As expected, both posterior means fit the data well, but the empirical Bayes estimate is more smooth.  That the empirical Bayes credible band is also a bit narrower than that of the Dirichlet process mixture should not be a concern based on the simulation results above.   
\end{ex}

\begin{figure}[t]
\begin{center}
\subfigure[Data and empirical Bayes]{\scalebox{0.6}{\includegraphics{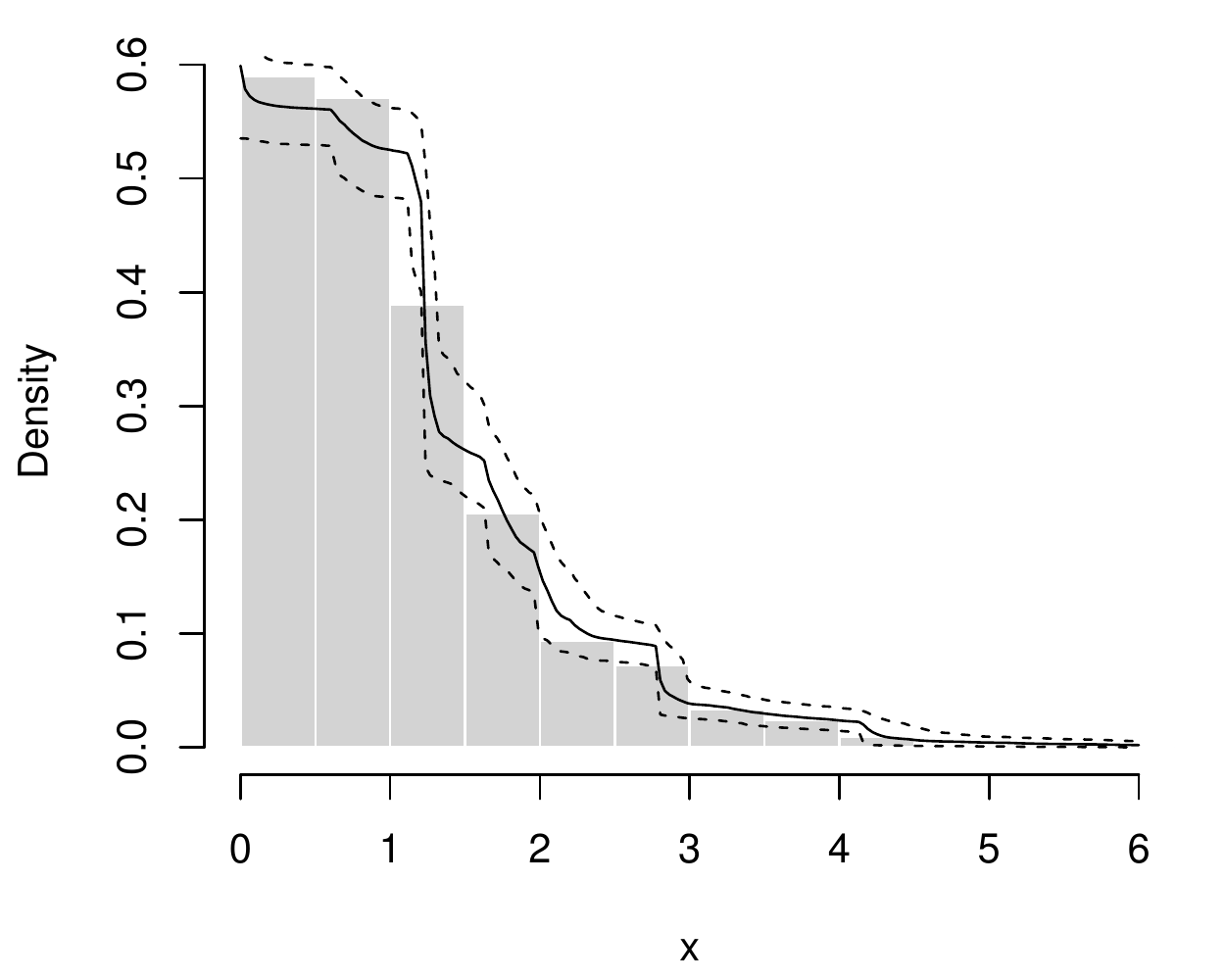}}}
\subfigure[Data and Dirichlet process mixture]{\scalebox{0.6}{\includegraphics{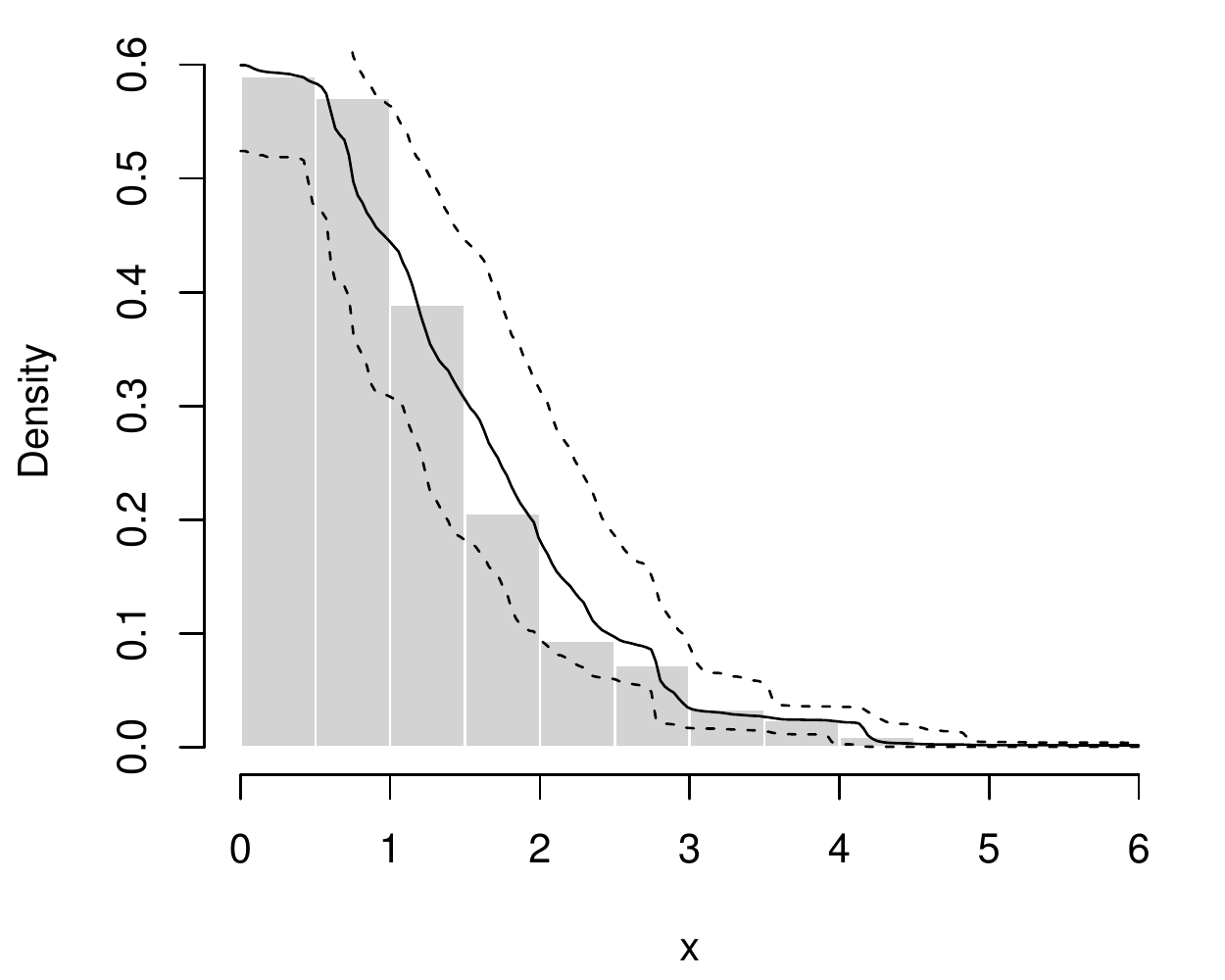}}}
\end{center}
\caption{Example~\ref{ex:norwayfire}, with $n=820$---same results as Figure~\ref{fig:silverman}.}
\label{fig:norwayfire}
\end{figure}

\begin{ex}
\label{ex:norwayfire}
The Norwegian fire claims data is a common example in the actuarial science literature \citep[e.g.,][]{brazauskas.kleefeld.2016}.  I consider $n=820$ fire loss claims exceeding $500$ thousand Norwegian krones during the year 1988.  Figure~\ref{fig:norwayfire} shows the data, the posterior mean densities, and the 95\% credible bands for the two methods.  In this case, the empirical Bayes posterior mean is a bit more rough than that of the Dirichlet process mixture, but arguably fits the data histogram better.  And the narrow credible band is expected since this data set is about ten times as large as that in Example~\ref{ex:silverman}.  
\end{ex}

% Example 2.3 in Groenboom's 2014 book
% Silverman 1986 has an example with data: 
%  https://ned.ipac.caltech.edu/level5/March02/Silverman/paper.pdf
% Turnbull and (Sujit) Ghosh looked at Silverman's data:
%  http://www.stat.ncsu.edu/information/library/papers/mimeo2650_Turnbull.pdf
% 1988 Norwegian fire loss data: log(loss / 500) must be positive, and shape looks monotone

\section{Conclusion}
\label{S:discuss}

This paper presents a unique approach to the specification of an empirical or data-dependent prior for nonparametric Bayesian-like inference on a density function.  The chief novelty is the centering of the prior on a suitable estimator, and this is particularly advantageous in the present context where the usual Kullback--Leibler property may fail.  The challenge is to choose the prior tails so that the corresponding posterior does not track the data too closely, and Theorems~\ref{thm:bounded}--\ref{thm:unbounded} provide sufficient conditions on the prior inputs to achieve the target posterior concentration rate.  Whether the proposed formulation can achieve valid uncertainty quantification, e.g., in the sense of \citet{szabo.vaart.zanten.2015}, remains an open question, but the numerical illustrations in Section~\ref{S:examples} are promising.  

%{\color{red} An important question being considered in the current Bayesian literature is calibration of posterior credible sets \citep[e.g.,][]{szabo.vaart.zanten.2015}.  A first step in proving this kind of strong calibration property is pinning down more precisely the constants involved in the rate derivations.  Here, and especially in the {\em clinical version} of the model (Remark~\ref{re:clinical}) considered in Section~\ref{S:examples}, the critical question is how to choose the proportionality constants attached to $c=c_n$ and $\delta=\delta_n$ in Theorems~\ref{thm:bounded}--\ref{thm:unbounded}.  An answer to this question would provide justification for the ``maximal spread'' claim in Remark~\ref{re:spread}, and this is a focus of some of my current work. }

\section*{Acknowledgment}

This work is partially supported by the National Science Foundation, DMS--1737933.  The author also thanks two anonymous the reviewers for their helpful feedback.

\bibliographystyle{apalike}
\bibliography{/Users/rgmarti3/Dropbox/Research/mybib}

\begin{thebibliography}{}

\bibitem[Balabdaoui and Wellner, 2007]{bala.wellner.2007}
Balabdaoui, F. and Wellner, J.~A. (2007).
\newblock Estimation of a {$k$}-monotone density: limit distribution theory and
  the spline connection.
\newblock {\em Ann. Statist.}, 35(6):2536--2564.

\bibitem[Birg\'e, 1989]{birge1989}
Birg\'e, L. (1989).
\newblock The {G}renander estimator: a nonasymptotic approach.
\newblock {\em Ann. Statist.}, 17(4):1532--1549.

\bibitem[Brazauskas and Kleefeld, 2016]{brazauskas.kleefeld.2016}
Brazauskas, V. and Kleefeld, A. (2016).
\newblock Modeling severity and measuring tail risk of {N}orwegian fire claims.
\newblock {\em N. Am. Actuar. J.}, 20(1):1--16.

\bibitem[Donnet et~al., 2018]{rousseau.etal.eb}
Donnet, S., Rivoirard, V., Rousseau, J., and Scricciolo, C. (2018).
\newblock Posterior concentration rates for empirical {B}ayes procedures with
  applications to {D}irichlet process mixtures.
\newblock {\em Bernoulli}, 24(1):231--256.

\bibitem[Ghosal et~al., 1999]{ggr1999}
Ghosal, S., Ghosh, J.~K., and Ramamoorthi, R.~V. (1999).
\newblock Posterior consistency of {D}irichlet mixtures in density estimation.
\newblock {\em Ann. Statist.}, 27(1):143--158.

\bibitem[Ghosal et~al., 2000]{ggv2000}
Ghosal, S., Ghosh, J.~K., and van~der Vaart, A.~W. (2000).
\newblock Convergence rates of posterior distributions.
\newblock {\em Ann. Statist.}, 28(2):500--531.

\bibitem[Ghosal and van~der Vaart, 2001]{ghosalvaart2001}
Ghosal, S. and van~der Vaart, A.~W. (2001).
\newblock Entropies and rates of convergence for maximum likelihood and {B}ayes
  estimation for mixtures of normal densities.
\newblock {\em Ann. Statist.}, 29(5):1233--1263.

\bibitem[Grenander, 1956]{grenander}
Grenander, U. (1956).
\newblock On the theory of mortality measurement. {II}.
\newblock {\em Skand. Aktuarietidskr.}, 39:125--153 (1957).

\bibitem[Groeneboom, 1985]{groeneboom}
Groeneboom, P. (1985).
\newblock Estimating a monotone density.
\newblock In {\em Proceedings of the {B}erkeley {C}onference in {H}onor of
  {J}erzy {N}eyman and {J}ack {K}iefer, {V}ol.\ {II} ({B}erkeley, {C}alif.,
  1983)}, Wadsworth Statist./Probab. Ser., pages 539--555, Belmont, CA.
  Wadsworth.

\bibitem[Groeneboom and Jongbloed, 2014]{groeneboom.jongbloed.book}
Groeneboom, P. and Jongbloed, G. (2014).
\newblock {\em Nonparametric Estimation under Shape Constraints}, volume~38 of
  {\em Cambridge Series in Statistical and Probabilistic Mathematics}.
\newblock Cambridge University Press, New York.

\bibitem[Kalli et~al., 2011]{kalli.griffin.walker.2011}
Kalli, M., Griffin, J.~E., and Walker, S.~G. (2011).
\newblock Slice sampling mixture models.
\newblock {\em Stat. Comput.}, 21(1):93--105.

\bibitem[Klaus and Strimmer, 2015]{fdrtool}
Klaus, B. and Strimmer, K. (2015).
\newblock {\em fdrtool: Estimation of (Local) False Discovery Rates and Higher
  Criticism}.
\newblock R package version 1.2.15.

\bibitem[Kosorok, 2008]{kosorok2008}
Kosorok, M.~R. (2008).
\newblock Bootstrapping in {G}renander estimator.
\newblock In {\em Beyond {P}arametrics in {I}nterdisciplinary {R}esearch:
  {F}estschrift in {H}onor of {P}rofessor {P}ranab {K}. {S}en}, volume~1 of
  {\em Inst. Math. Stat. (IMS) Collect.}, pages 282--292. Inst. Math. Statist.,
  Beachwood, OH.

\bibitem[Martin et~al., 2017]{martin.mess.walker.eb}
Martin, R., Mess, R., and Walker, S.~G. (2017).
\newblock Empirical {B}ayes posterior concentration in sparse high-dimensional
  linear models.
\newblock {\em Bernoulli}, 23(3):1822--1847.

\bibitem[Martin and Walker, 2014]{martin.walker.eb}
Martin, R. and Walker, S.~G. (2014).
\newblock Asymptotically minimax empirical {B}ayes estimation of a sparse
  normal mean vector.
\newblock {\em Electron. J. Stat.}, 8(2):2188--2206.

\bibitem[Martin and Walker, 2017]{martin.walker.deb}
Martin, R. and Walker, S.~G. (2017).
\newblock Empirical priors for target posterior concentration rates.
\newblock Unpublished manuscript, {\tt arXiv:1604.05734}.

\bibitem[Prakasa~Rao, 1969]{prakasarao}
Prakasa~Rao, B. L.~S. (1969).
\newblock Estimation of a unimodal density.
\newblock {\em Sankhy\=a Ser. A}, 31:23--36.

\bibitem[Robertson et~al., 1988]{robertson.wright.dykstra.1988}
Robertson, T., Wright, F.~T., and Dykstra, R.~L. (1988).
\newblock {\em Order {R}estricted {S}tatistical {I}nference}.
\newblock Wiley Series in Probability and Mathematical Statistics: Probability
  and Mathematical Statistics. John Wiley \& Sons, Ltd., Chichester.

\bibitem[Rousseau and Szabo, 2017]{rousseau.szabo.2017}
Rousseau, J. and Szabo, B. (2017).
\newblock Asymptotic behaviour of the empirical {B}ayes posteriors associated
  to maximum marginal likelihood estimator.
\newblock {\em Ann. Statist.}, 45(2):833--865.

\bibitem[Salomond, 2014]{salomond2014}
Salomond, J.-B. (2014).
\newblock Concentration rate and consistency of the posterior distribution for
  selected priors under monotonicity constraints.
\newblock {\em Electron. J. Stat.}, 8(1):1380--1404.

\bibitem[Schwartz, 1965]{schwartz1965}
Schwartz, L. (1965).
\newblock On {B}ayes procedures.
\newblock {\em Z. Wahrs. verw. Geb.}, 4:10--26.

\bibitem[Sen et~al., 2010]{sen.banerjee.woodroofe.2010}
Sen, B., Banerjee, M., and Woodroofe, M. (2010).
\newblock Inconsistency of bootstrap: the {G}renander estimator.
\newblock {\em Ann. Statist.}, 38(4):1953--1977.

\bibitem[Silverman, 1986]{silverman}
Silverman, B.~W. (1986).
\newblock {\em Density Estimation for Statistics and Data Analysis}.
\newblock Chapman \& Hall, London.

\bibitem[Szab\'o et~al., 2015]{szabo.vaart.zanten.2015}
Szab\'o, B., van~der Vaart, A.~W., and van Zanten, J.~H. (2015).
\newblock Frequentist coverage of adaptive nonparametric {B}ayesian credible
  sets.
\newblock {\em Ann. Statist.}, 43(4):1391--1428.

\bibitem[Szab{\'o} et~al., 2013]{szabo.vaart.zanten.2013}
Szab{\'o}, B.~T., van~der Vaart, A.~W., and van Zanten, J.~H. (2013).
\newblock Empirical {B}ayes scaling of {G}aussian priors in the white noise
  model.
\newblock {\em Electron. J. Stat.}, 7:991--1018.

\bibitem[van~der Vaart and Wellner, 1996]{vaartwellner1996}
van~der Vaart, A.~W. and Wellner, J.~A. (1996).
\newblock {\em Weak Convergence and Empirical Processes}.
\newblock Springer-Verlag, New York.

\bibitem[Walker, 2007]{walker2007.slice}
Walker, S.~G. (2007).
\newblock Sampling the {D}irichlet mixture model with slices.
\newblock {\em Comm. Statist. Simulation Comput.}, 36(1-3):45--54.

\bibitem[Walker et~al., 2007]{walker2007}
Walker, S.~G., Lijoi, A., and Pr{\"u}nster, I. (2007).
\newblock On rates of convergence for posterior distributions in
  infinite-dimensional models.
\newblock {\em Ann. Statist.}, 35(2):738--746.

\bibitem[Williamson, 1956]{williamson1956}
Williamson, R.~E. (1956).
\newblock Multiply monotone functions and their {L}aplace transforms.
\newblock {\em Duke Math. J.}, 23:189--207.

\bibitem[Woodroofe and Sun, 1993]{woodroofesun}
Woodroofe, M. and Sun, J. (1993).
\newblock A penalized maximum likelihood estimate of {$f(0+)$} when {$f$} is
  nonincreasing.
\newblock {\em Statist. Sinica}, 3(2):501--515.

\bibitem[Wu and Ghosal, 2008]{wu.ghosal.2008}
Wu, Y. and Ghosal, S. (2008).
\newblock Kullback {L}eibler property of kernel mixture priors in {B}ayesian
  density estimation.
\newblock {\em Electron. J. Stat.}, 2:298--331.

\end{thebibliography}

\end{document}